\documentclass[10pt]{article}
\usepackage{amsmath, amssymb,amsthm}
\newtheorem{theo}{{\bf{Theorem}}}[section]
\newtheorem{prop}[theo]{{\bf Proposition}}
\newtheorem{lem}[theo]{{\bf Lemma}}

\newtheorem{rem}{{\bf Remark}}[section]

\newtheorem{remark}[rem]{Remark}

\renewcommand{\proof}{\noindent{\bf Proof.\ }}
\newcommand{\no}{\nonumber}
\newcommand{\noi}{\noindent}
\newcommand{\txt}{\textrm}
\newcommand{\pa}{\partial}

\newcommand{\RR}{\mathbb{R}}
\newcommand{\si}{\sigma}

\newcommand{\la}{\lambda}
\newcommand{\La}{\Lambda}
\newcommand{\calF}{\mathcal{F}}

\newcommand{\sig}{\sigma}
\newcommand{\tab}{\hspace*{0.3in}}

\newcommand{\vf}{\varphi}

\catcode`\@=11
\catcode`\@=12

\begin{document}

\title{A system of non-local parabolic PDE and application to option pricing}
\author{Anindya Goswami\thanks{IISER, Pune 411008, India; email: anindya@iiserpune.ac.in}\qquad\qquad Jeeten Patel
\thanks{IISER, Pune 411008 India; email: jeetenp03@gmail.com}\qquad\qquad Poorva Shevgaonkar \thanks{IIT Kharagpur, India; email: poorvashevgaonkar@gmail.com}\\ }


\maketitle

{\bf Abstract:} This paper includes a proof of well-posedness of an initial-boundary value problem involving a system of non-local parabolic partial differential equation(PDE) which naturally arises in the study of derivative pricing in a generalized market model which is known as a semi-Markov modulated geometric Brownian motion(GBM) model. We study the well-posedness of the problem via a Volterra integral equation of second kind. A probabilistic approach, in particular the method of conditioning on stopping times is used for showing the uniqueness.


{\bf Classification No.:} 60K15, 91B30, 91G20, 91G60.
\section{Introduction}
In the Black-Scoles-Merton model of option pricing the dynamics of stock price $\{S_t\}_{t\ge 0}$ is assumed to be given by a
geometric Brownian motion (GBM), i.e.,
$$ dS_t =  \mu S_t dt +\sigma S_t dW_t, \tab S_0>0
$$
where $\mu$ and $\sigma$ are positive constants and $ \{W_t\}_{t\geq 0}$ is a standard Wiener process. Existence of a deterministic
 risk free asset with constant growth rate $r$ is also assumed. Under these model assumptions, the European call option price
 function satisfies a parabolic partial differential equation
\begin{eqnarray*}
\no \frac{\partial}{\partial t} \varphi(t, s ) + r  s \frac{\partial} {\partial s} \varphi(t, s ) + \frac{1}{2} \si^2
 s^2 \frac{\partial^2} {\partial s^2} \varphi(t, s )
 = r  ~ \varphi(t, s ).
\end{eqnarray*}
This is known as Black-Scholes-Merton partial differential equation (B-S-M PDE), whose solution can be obtained explicitly \cite{BS}. However, empirical studies of financial assets suggest that the assumption of constant or deterministic $\sigma$, the volatility coefficient, is rather unrealistic. To overcome this drawback, there has been an increasing interest in market models where $\sigma$ is taken as a Markov process. In stochastic volatility models, the square of $\sigma$ is taken as an Ito diffusion \cite{HE}, whereas in regime switching models $\sigma$ is driven by a finite-state continuous-time Markov chain \cite{BAS}, \cite{BE}, \cite{DES}, \cite{DKR}, \cite{ECS}, \cite{GoS}, \cite{GZ}, \cite{RR} and \cite{MR}. The resulting market models for both of theses types are incomplete. Option pricing in such markets has been studied by several authors using different approaches. Using locally risk minimizing approach, the European call option price in a regime switching market is shown to satisfy a generalized B-S-M PDE \cite{BAS}, \cite{DES}, \cite{DKR}.

\noi Locally risk minimizing option pricing in a regime switching model with semi-Markov regimes is studied in \cite{AGMKG}. There the market parameters, $\mu$, $\sigma$ and $r$ are driven by a finite-state semi-Markov process $\{X_t\}_{t\ge 0}$. Sojourn or holding times in each state in a finite-state continuous-time Markov chain are distributed as exponential random variables, whereas that could be any positive random variable for semi-Markov case. Thus the class of semi-Markov processes subsumes the class of Markov chains. There are some statistical results in the literature (see \cite{JH} and the references therein) which emphasize the advantage of use of semi-Markov switching models over simple homogeneous Markov switching models. For example it is mainly useful to deal with the impact of a changing environment (i.e. the business cycle), which exhibits duration dependence. Suppose $\{X_t\}$ is a semi-Markov process and the stock price $\{S_t\}$ is given by
$$ dS_t =  \mu (X_t) S_t dt +\sigma (X_t) S_t dW_t, \tab S_0>0.
$$
Then  it is shown in \cite{AGMKG} that the call option price function satisfies a non local system of parabolic PDE
\begin{eqnarray}
\no \frac{\partial}{\partial t} \varphi(t, s, i, y)+ \frac{\partial} {\partial y} \varphi(t, s, i, y) + r( i ) s
\frac{\partial} {\partial s} \varphi(t, s, i, y) + \frac{1}{2} \si^2( i ) s^2 \frac{\partial^2} {\partial s^2} \varphi(t, s, i, y)\\
\no +\frac{f(y\mid i)}{1-F(y\mid i)}\sum_{j\neq i}p_{ij}[\varphi(t, s, j, 0) -\varphi(t, s, i, y)] = r( i ) ~ \varphi(t, s, i, y),
\end{eqnarray}
\noi defined on
\begin{equation}\label{D}
\mathcal{D}:= \{ (t, s, i, y)\in (0,T)\times \mathbb{R_+}\times \mathcal{X}\times (0,T) \mid y \in (0,t)\},
\end{equation}
\noi and with conditions
\begin{eqnarray}\label{18}
\no \lim_{s \downarrow 0} \varphi (t, s, i, y)&=& 0 \tab \forall t \in [0,T]\\
\varphi (T, s, i, y)&=& (s-K)^+; \tab s\in (0,\infty); \tab 0 \le y\le T ; \tab i= 1, 2, \ldots, k
\end{eqnarray}
where $k$ is the number of all possible states of $X_t$, $r(X_t)$ is the spot rate of interest at time $t$, $(p_{ij})$ are the transition probabilities to state $j$ form $i$, $F(\cdot\mid i)$ is the conditional distribution function of holding time, $f(y\mid i)=\frac{d}{dy}F(y\mid i)$, $K$ is the strike price and $T$ is the maturity time.

\noi In this paper we replace the semi-Markov process $\{X_t\}$ by a more general  class of age-dependent processes which is a much wider class of switching than that appear in \cite{AGMKG}. An age-dependent process $\{X_t\}_{t\ge 0}$ is specified by its instantaneous transition rate which is a collection of measurable functions $\la_{ij}:[0,\infty)\to (0,\infty)$ where $i\neq j \in \mathcal{X}:=\{1,2,\ldots,k\}$. Indeed $\{X_t\}_{t\geq0}$ is defined as the strong solution of the following system of stochastic integral equations
\begin{eqnarray}\label{1}
  X_t &=& X_0 + \int_{(0,t]}\int_{\mathbb{R}} h(X_{u-}, Y_{u-},z)\wp(du,dz)\\
  \no Y_t &=& t- \int_{(0,t]} \int_{\mathbb{R}} g(X_{u-}, Y_{u-},z) \wp(du,dz)
\end{eqnarray}
where $\wp (dt,dz)$ is the Poisson random measure with intensity $dtdz$, independent
of $X_0$ and $$  h(i,y,z) := \sum_{j \in \mathcal{X}, j\neq i} (j-i) 1_{\Lambda_{ij}(y)}(z),
 ~~~  g(i,y,z) := \sum_{j \in \mathcal{X}, j\neq i} y 1_{\Lambda_{ij}(y)}(z)$$
where  $\Lambda_{ij}(y)$ are the consecutive (with respect to the lexicographic ordering on $\mathcal{X}\times \mathcal{X}$) left closed and right open intervals of the real line, each having length $\lambda_{ij} (y)$. We refer \cite{GhS} for more details about this kind of pure jump processes. Under some smoothness and tail assumptions on $\lambda_{ij} (y)$ and independence of $W$ and $\wp$, we obtain the following equation of locally risk minimizing price of call option
\begin{eqnarray}\label{p1}
\no \frac{\partial}{\partial t} \varphi(t, s, i, y)+ \frac{\partial} {\partial y}
\varphi(t, s, i, y) + r( i ) s \frac{\partial} {\partial s} \varphi(t, s, i, y)
+ \frac{1}{2} \si^2( i ) s^2 \frac{\partial^2} {\partial s^2} \varphi(t, s, i, y)\\
+\sum_{j\neq i}\la_{ij}(y)\big(\varphi(t, s, j, 0) -\varphi(t, s, i, y)\big) = r( i ) ~ \varphi(t, s, i, y),
\end{eqnarray}
\noi defined on $\mathcal{D}$ as in \eqref{D} and with conditions as in \eqref{18}. Some of the special cases of this equation appear in \cite{BAS}, \cite{DES}, \cite{RR}, \cite{AGMKG}, \cite{MR}, and \cite{DKR} for pricing a European contingent claim under certain regime switching market assumptions. Owing to the simplicity of the special cases, generally authors refer to some standard results in the theory of parabolic PDE for existence and uniqueness issues. But in its general form which arises in this paper, no such ready reference is available. In this connection we would like to emphasize that \eqref{p1} is a non-local degenerate parabolic PDE on a non-rectangular domain. Therefore, we produce a self contained proof using Banach fixed point theorem. These we accomplish in two steps. First we consider a Volterra integral equation(IE) of second kind and establish existence and uniqueness result of that. Then we show in a couple of propositions, that the PDE and the integral equation(IE) problems are ``equivalent". Thus we obtain the existence and uniqueness of the PDE in Theorem \ref{theo1}. Some further properties, viz the positivity and growth property are also obtained. It is also shown here that the partial derivative of the solution constitutes the optimal hedging strategy of the corresponding claim. We further show that the partial derivative of $\varphi$, can be written as an integration involving $\varphi$ which enables one to develop a robust numerical scheme to compute the Greeks. This study paves the way for addressing many other interesting problems involving this new set of PDEs.

\noi The rest of this paper is arranged in the following manner. The market model assumption and the pricing approach is described in Section 2. In Sections 3 and 4 the well-posedness of the Cauchy problem \eqref{p1}is established. In Section 5, the well-posedness result is used to solve the pricing and hedging problem. We conclude our paper in Section 6 with a few remarks.
\section{Market model}
\noi We consider the maps $\la_{ij}:[0,\infty)\to (0,\infty)$ where $i\neq j \in \mathcal{X}:=\{1,2,\ldots,k\}$ and define $ \la_{ii}(y):=-\sum_{j\neq i \\ j\in \mathcal{X}}\la_{ij}(y)$ for all $i\in \mathcal{X}$ and $y \in [0,\infty)$. We further assume the following.
\begin{itemize}
\item[(A1)](i) For each $i$ and $j$, $\la_{ij}$ is a differentiable function, \\
(ii)$\lim_{y\to \infty}\La_i(y) =\infty$, where $\La_i(y):=\int_0^y|\la_{ii}(y)|dy$.
\end{itemize}
Define $F(y\mid i):= 1-\exp(-\La_i(y))$, $f(y\mid i):= \frac{d}{dy}F(y\mid i)$ and for each $i\neq j$, $p_{ij}(y):=\frac{\la_{ij}(y)}{|\la_{ii}(y)|}$ with $p_{ii}(y)=0$ for all $i$ and $y$. Set $\hat{p}_{ij}:=\int_0^\infty p_{ij}(y) dF(y\mid i)$.
\begin{itemize}
\item[(A1)](iii) The matrix $(\hat{p}_{ij})$ is irreducible.
\end{itemize}
Let $(\Omega,\mathcal{F},P)$ be the underlying complete probability space where a Poisson random measure $\wp$ and a standard Wiener process $W$ are defined and they are independent. Let $(X_t, Y_t)$ be the strong solution of system of equations \eqref{1} where the maps $\la_{ij}$ satisfy (A1) (i)-(iii). From \eqref{1} it is apparent that $X_t$ is a right continuous (since the integrations are over $(0,t]$) jump process having left limits and taking values in $\mathcal{X}$. Let $T_n$ denote the time of $n$th jump of $X_t$, whereas $T_0:= 0$ and $\tau_n:=T_n-T_{n-1}$. For a fixed $t$, let $n(t):= \max\{n: T_n \leq t \}$. Thus $T_{n(t)}\leq t < T_{n(t)+1}$ and $Y_t = t-T_{n(t)}$. It is shown in \cite{GhS} that $F(y\mid i)$ is the  cumulative distribution function of holding time and $p_{ij}(y)$ is the conditional probability that $X$ transits to $j$ given the fact that it is at $i$ for a duration of $y$. From (A1)(ii), $\lim_{y\to \infty}F(y\mid i)=1$. Hence, sojourn times are finite almost surely. We also note that for $i\neq j$, $\la_{ij}(y)=p_{ij}(y)\frac{f(y\mid i)}{1-F(y\mid i)}$ hold. Thus the instantaneous transition rate is given by the maps $\la_{ij}$.

\noi Let $\{B_t\}_{t\geq0}$ be the price of money market account at time $t$ where, spot interest rate is $r_t=r(X_t)$ and $B_0=1$. Thus we have $B_t = e^{\int_0^t r(X_{u}) du}$. Let $\{S_t\}_{t\geq 0}$ be the price process of the stock, which is governed by a semi-Markov modulated GBM i.e.,
\begin{equation}\label{8}
dS_t = \mu(X_{t})S_tdt +\sigma(X_{t})S_t dW_t,\tab S_0>0
\end{equation}
where $\mu : \mathcal{X} \to \mathbb{R}$ is the drift coefficient and $\sigma : \mathcal{X} \to (0, \infty)$ corresponds to the volatility.
Let $\mathcal{F}_t$ be a filtration of $\mathcal{F}$ satisfying usual hypothesis and right continuous version of the filtration generated by $X_t$ and $S_t$. Clearly the solution of the above SDE is an $\mathcal{F}_t$ semimartingale with almost sure continuous paths. We assume that $X_t$ is observed.

\noi An \emph{admissible strategy} is a dynamic allocation to these assets and is defined as a predictable process $\pi=\{\pi_t=(\xi_t,\varepsilon_t), 0\leq t\leq T\}$ which satisfies conditions, given in $(A2)$ below.
\begin{itemize}
\item[(A2)](i) $\xi_t$ is square integrable w.r.t $S_t$,\\
(ii) $E(\varepsilon^2_t)<\infty$,\\
(iii) $\exists a>0$ s.t. $P(V_t\geq-a,t\in[0,T])=1$ where $V_t = \xi_tS_t + \varepsilon_tB_t$.
\end{itemize}
It is shown in \cite{FS} that if the market is arbitrage free, the existence of an optimal strategy for hedging an $\calF_T$ measurable claim $H$, is equivalent to the existence of F\"{o}llmer-Schweizer decomposition of discounted claim $H^*:= B^{-1}_T H$ in the form
\begin{equation}\label{eq1}
H^*=H_0+\int^{T}_{0}{\xi^{H^*}_t}dS^*_t+L^{H^*}_T
\end{equation}
where $H_0\in L^2(\Omega,\mathcal{F}_0,P), L^{H^*}=\{L^{H^*}_t\}_{0\leq t\leq T}$ is a square integrable martingale starting with zero and orthogonal to the martingale part of $S_t$, and $\xi^{H^*}=\{\xi^{H^*}_t\}$ satisfies A2 (i). Further $\xi^{H^*}$ appeared in the decomposition, constitutes the optimal strategy. Indeed the optimal strategy $\pi=(\xi_t,\varepsilon_t)$ is given by
\begin{eqnarray}\label{11a}
 \nonumber \xi_t &:=& \xi^{H^*}_t,\\
V^*_t &:=& H_0+\int^{t}_{0}{\xi_u}dS^*_u+L^{H^*}_t,\\
 \nonumber \varepsilon_t &:=& V^*_t-\xi_tS^*_t,
\end{eqnarray}
\noi and $B_t V^*_t$ represents the \emph{locally risk minimizing price} at time $t$ of the claim $H$. Hence the F\"{o}llmer-Schweizer decomposition is the key thing to verify for settling the pricing and hedging problems in a given market (in particular when it is incomplete). We refer to \cite{S2} for more details.

\noi It can be shown in the similar line of \cite{AGMKG} that this market model admits existence of an equivalent martingale measure. Hence under admissible strategies the market is arbitrage free. To price a claim $H$ of European type in the above incomplete market, we would consider the locally risk minimizing pricing approach by F\"{o}llmer and Schweizer, i.e., decomposition of type \eqref{11a} and then show that the strategy, thus obtained is admissible.
\section{Existence}
Consider the following initial boundary value problem which is known as B-S-M PDE for each $i$
\begin{equation}\label{eq2}
\frac{\partial\eta_{i}(t,s)}{\partial t} + r(i)s\frac{\partial\eta_{i}(t,s)}{\partial s}+\frac{1}{2}\sigma^2(i)s^2 \frac{\partial^2\eta_{i}(t,s)}{\partial s^2} = r(i)\eta_{i}(t,s)
\end{equation}
for $(t,s)\in (0,T)\times (0,\infty)$ and $\eta_{i}(T,s)=(s-K)^+$,  $\lim_{s \downarrow 0} \eta_{i}(t, s)= 0$ for all $t\in [0,T]$. This has a classical solution with at most linear growth (see \cite{KK}). We also introduce a log normal probability density function
\begin{equation*}
\alpha(x;s,i,v):=\frac{e^{-\frac{1}{2}\beta^2}}{\sqrt{2\pi}x\sigma(i)\sqrt{v}},~~~~
\beta(x,s,i,v):=\frac{\ln\left(\frac{x}{s}\right)-\left(r(i)-\frac{\si^2(i)}{2}\right)v}{\sigma(i)\sqrt{v}}.
\end{equation*}
By a direct calculation one has
\begin{equation}\label{10ae}
\beta\frac{\partial \beta}{\partial v}+r(i)\frac{\beta}{\sigma(i) \sqrt{v}}+\frac{1}{2}\frac{\beta^2}{v}-\frac{\sigma (i)\beta}{2\sqrt{v}}=0.
\end{equation}
Set $B:=\left\{\varphi:\mathcal{D}\rightarrow \RR, \mathrm{continuous} \mid \lim_{s \downarrow 0} \varphi(\cdot,s,\cdot,\cdot)=0, ~\|\varphi\|:=\sup_{\mathcal{D}}\mid\frac{\varphi(t,s,i,y)}{1+s}\mid < \infty\right\}.
$
\begin{lem}\label{lm2}
Consider the following integral equation
\begin{eqnarray}
\no\varphi(t,s,i,y)&=&\frac{1- F(T-t+y\mid i)}{1-F(y\mid i)} \eta_{i}(t,s)+\int_0^{T-t} e^{-r(i)v} \frac{f(y+v\mid i)} {1-F(y\mid i)} \times\\
\label{34} && \sum_{j\neq i} p_{ij}(y+v) \int_0^{\infty} \varphi(t+v,x,j,0)\alpha(x;s,i,v) dx dv \\
\label{35} \txt{ with } \lim_{s \downarrow 0} \varphi(t,s,i,y)&=&0 ~\forall t \in [0,T], i\in \chi,~y\in [0,t]
\end{eqnarray}
where $ \eta$ is as in \eqref{eq2}. Then (i) the problem \eqref{34}-\eqref{35} has unique solution in $B$, and (ii) the solution of the integral equation is in $C^{1,2,1}(\mathcal{D})$, (iii) $\varphi(t,s,i,y)$ is nonnegative.
\end{lem}
\proof (i) We first note that a solution of \eqref{34}-\eqref{35} is a fixed point of the operator $A$ and vice versa, where
\begin{eqnarray}
\no A \varphi(t,s,i,y)&:=& \frac{1- F(T-t+y\mid i)}{1-F(y\mid i)}\eta_{i}(t,s)+\int_0^{T-t}e^{-r(i) v} \frac{f(y+v\mid i)}{1-F(y\mid i)}
 \sum_{j\neq i} p_{ij}(y+v)\\
 \no && \int^\infty _0 \varphi(t+v,x,j,0) \alpha(x;s,i,v)dx dv,
\end{eqnarray}

\noi Clearly $B$ is a closed subspace of a Banach space $(\mathcal{B}, \|\cdot \|)$, where $\mathcal{B}$ is the set of all continuous functions with at most linear growth in $s$ variable. It is also easy to check that for each $\vf\in B$, $A\vf: \mathcal{D}\to (0,\infty)$ is continuous and $\lim_{s \downarrow 0} A\varphi(\cdot,s,\cdot,\cdot)=0$. Now in order to show existence and uniqueness in the prescribed class, it is sufficient to show that $A$ is a contraction in $\mathcal{B}$. Because, then $A:B\to B$ is also a contraction and Banach fixed point theorem ensures existence and uniqueness of the fixed point in $B$. To this end, we need to show, for $\vf_1,\vf_2\in \mathcal{B}$, $||A\varphi_1-A\varphi_2|| \leq J||\varphi_1-\varphi_2||$ where $J<1$. Indeed
\begin{eqnarray*}
\|A\varphi_1-A\varphi_2\|&=&\sup_{\mathcal{D}}\bigg|\frac{ A\varphi_1-A\varphi_2}{1+s}\bigg|\\
&=&\sup_{\mathcal{D}}\bigg|\int^{T-t}_0 e^{-r(i) v} \frac{f(y+v\mid i)}{1-F(y\mid i)} \sum_{j\neq i} p_{ij}(y+v) \int^\infty_0 (\varphi_1-\varphi_2)(t+v,x,j,0)\frac{\alpha(x;s,i,v)}{1+s}dx dv\bigg|\\
&\leq& \sup_{\mathcal{D}} \bigg| \int^{T-t}_0 e^{-r(i) v} \frac{f(y+v\mid i)}{1-F(y\mid i)}\sum_{j\neq i} p_{ij}(y+v) \int^\infty_0 (1+x)~\sup_{\mathcal{D}}\bigg|\frac{\varphi_1-\varphi_2}{1+x}\bigg| \frac{\alpha(x;s,i,v)}{1+s}dx dv\bigg|\\
&=& \sup_\mathcal{D} \bigg| \int^{T-t}_0 e^{-r(i) v} \frac{f(y+v \mid i)}{1-F(y \mid i)}\|\varphi_1- \varphi_2 \| \frac{a(s)}{1+s}dv\bigg|
\end{eqnarray*}
where, $a(s)= \int^\infty_0 (1+x) \alpha(x;s,i,v) dx = 1 + e^{\ln s + \left( r(i)-\frac{\sigma^2(i)}{2}\right)v+\frac{\sigma^2(i) v}{2}}= 1+ s e^{r(i)v}.$

\noi Thus, $\|A\varphi_1-A\varphi_2\| \, \leq\, J\|\varphi_1-\varphi_2\|$ where,
\begin{eqnarray*}
J&=&\sup_\mathcal{D} \bigg|\int^{T-t}_0e^{-r(i) v}\frac{f(y+v \mid i)}{1-F(y \mid i)}\frac{1+se^{r(i)v}}{1+s}dv\bigg|\\
&\leq &\sup_\mathcal{D}\bigg( \frac{1}{1-F(y\mid i)}\int^{T-t}_0 f(y+v|i)dv\bigg)\\
&=&\sup_\mathcal{D}\bigg(\frac{F\left( y+T-t\mid i\right) -F(y|i)}{1-F\left(y|i\right)}\bigg)\\
&<&\frac{1-F(y|i)}{1-F(y|i)}~=~1
\end{eqnarray*}
using $r(i)\geq 0$ and (A1).

\noi (ii) Using (A1) and smoothness of $\eta_i$ for each $i$, the first term on the right hand side is in $C^{1,2,1}(\mathcal{D})$. Under assumption (A1) the fact, the second term is continuous differentiable in $y$ and twice continuously differentiable in s, follows immediately. The continuous differentiability in $t$ follows from the fact that the term $\varphi(t+v,x,j,0)$ is multiplied by $C^1((0,\infty))$ functions in $v$ and then integrated over $v\in(0,T-t)$. Hence $\varphi(t,s,i,y)$ is in $C^{1,2,1}(\mathcal{D})$.

\noi (iii)  We note that the problem \eqref{eq2} has a closed form solution and that is nonnegative. Hence nonnegativity of $\vf$ follows from nonnegativity of all the coefficients of Volterra equation \eqref{34}.\qed

\begin{prop}\label{theo3} The unique solution of \eqref{34}-\eqref{35} also solves the initial boundary value problem \eqref{p1}-\eqref{18}.
\end{prop}
\proof
Let $\varphi$ be the solutions of \eqref{34}-\eqref{35}. Thus using \eqref{34}, $\varphi(T,s,i,y)=\eta_{i}(T,s)=(s-K)^+$, i.e., the condition \eqref{18} holds. From Lemma \ref{lm2} (ii), $\varphi$ is in $C^{1,2,1}(\mathcal{D})$. Hence we can perform the partial differentiations w.r.t. $t$ and $y$ on the both sides of \eqref{34}. We obtain
\begin{eqnarray}\label{2ae}
\no \frac{\partial}{\partial t} \varphi(t, s, i, y)&=&\frac{f(T-t+y|i)}{1-F(y\mid i)}\eta_{i}(t,s)+\frac{ 1- F(T-t+y\mid i)}{\left(1-F(y\mid i)\right)}\frac{\partial\eta_{i}(t,s)}{\partial t}- e^{-r(i)(T-t)}\frac{f(y+T-t\mid i)}{1-F(y\mid i)}\\
&&\no\sum_{j\neq i} p_{ij}(y+T-t)\int_0^\infty \varphi(T,x,j,0)\alpha(x;s,i,T-t)dx+\int^{T-t}_0 e^{-r(i)v}\frac{f(y+v\mid i)}{1-F(y\mid i)} \\
&&\sum_{j\neq i} p_{ij}(y+v)\int^\infty_0 \frac{\partial \varphi}{\partial t}(t+v,x,j,0) \alpha(x;s,i,v) dx dv
\end{eqnarray}
by differentiating w.r.t. $t$ under the sign of integral. Now, after taking partial derivative w.r.t. $y$ on both sides of \eqref{34}, we get 
\begin{eqnarray*}
\no \frac{\partial}{\partial y} \varphi(t, s, i, y)&=&\no -\frac{f(T-t+y \mid i)}{1-F(y\mid i)}\eta_{i}(t,s)+\frac{1- F(T-t+y\mid i)}{\left(1-F(y\mid i)\right)^2}f(y|i)\eta_{i}(t,s)+\frac{f(y|i)}{1-F(y\mid i)}\times\\
&&\no \Big(\varphi(t,s,i,y)- \frac{1- F(T-t+y \mid i)}{1-F(y\mid i)}\eta_{i}(t,s)\Big)\\
&& +\sum_{j\neq i} \int_0^{T-t} e^{-r(i)v} \frac{{q_{ij}}'(y+v)} {1-F(y\mid i)} \int_0^{\infty} \varphi(t+v,x,j,0)\alpha(x;s,i,v) dx dv.
\end{eqnarray*}
where $q'(u):= \frac{d}{du}q(u)$ and $q_{ij}(u)=f(u| i)p_{ij}(u)$. For further simplification, we would simplify the last additive term in the right side using integration by parts w.r.t. $v$ where $q'$ is treated as second function. We would also use the following identities $\frac{\pa}{\pa t}\varphi(t+v,x,j,0)= \frac{\pa}{\pa v}\varphi(t+v,x,j,0)$ and $\frac{\partial \alpha}{\partial v}= -\alpha\left(\beta \frac{\partial \beta}{\partial v}+\frac{1}{2v}\right)$ to obtain
\begin{eqnarray}\label{3ae}
\no \frac{\partial}{\partial y} \varphi(t, s, i, y)&=&\no -\frac{f(T-t+y \mid i)}{1-F(y\mid i)}\eta_{i}(t,s)+\frac{1- F(T-t+y\mid i)}{\left(1-F(y\mid i)\right)^2}f(y|i)\eta_{i}(t,s)+\frac{f(y|i)}{1-F(y\mid i)}\times\\
&&\no \Big(\varphi(t,s,i,y)- \frac{1- F(T-t+y \mid i)}{1-F(y\mid i)}\eta_{i}(t,s)\Big)+e^{-r(i)(T-t)}\frac{f(T-t+y\mid i)}{1-F(y\mid i)} \\
&&\no\sum_{j\neq i} p_{ij}(y+T-t)\int^\infty_0 \varphi(T,x,j,0)\alpha(x;s,i,T-t)dx-\frac{f(y\mid i)}{1-F(y\mid i)}\sum_{j\neq i} p_{ij}(y)\varphi(t,s,j,0)\\
&&\no-\int^{T-t}_0 e^{-r(i) v} \frac{f(y+v\mid i)}{1-F(y\mid i)}\int^\infty_0 \alpha(x;s,i,v) \Bigg\{-r(i) \sum_{j\neq i} p_{ij}(y+v) \varphi(t+v,x,j,0)\\
&&-\sum_{j\neq i} p_{ij}(y+v) \varphi(t+v,x,j,0)\left(\beta \frac{\partial \beta}{\partial v}+\frac{1}{2v}\right)+\sum_{j\neq i} p_{ij}(y+v) \frac{\partial{\varphi(t+v,x,j,0)}}{\partial t}\Bigg\}dx dv.
\end{eqnarray}
\noi By adding equations \eqref{2ae} and \eqref{3ae}, we get
\begin{eqnarray}\label{6ae}
\no\lefteqn{\frac{\partial}{\partial t} \varphi(t, s, i, y)+\frac{\partial} {\partial y} \varphi(t, s, i, y)}\\
&=&\no\frac{ 1- F(T-t+y\mid i)}{1-F(y\mid i)}\frac{\partial\eta_{i}(t,s)}{\partial t}+ \frac{f(y|i)}{1-F(y\mid i)}\left(
\varphi(t,s,i,y)-\sum_{j\neq i} p_{ij}(y)\varphi(t,s,j,0)\right) + \int^{T-t}_0 e^{-r(i)v}\\
&& \frac{f(y+v|i)}{1-F(y|i)}\sum_{j\neq i} p_{ij}(y+v) \int_0^\infty \varphi(t+v,x,j,0)\alpha(x;s,i,v)\left(r(i)+\beta\frac{\partial \beta}{\partial v}+\frac{1}{2v}\right)dv dx.
\end{eqnarray}
\noi Now we differentiate both sides of \eqref{34} w.r.t. $s$ once and twice respectively and obtain
\begin{eqnarray}\label{4ae}
\no \frac{\partial}{\partial s} \varphi(t, s, i, y)&=& \frac{ 1- F(T-t+y\mid i)}{1-F(y\mid i)}\frac{\partial\eta_{i}(t,s)}{\partial s}+\int_0^{T-t} e^{-r(i)v}\frac{f(y+v\mid i)} {1-F(y\mid i)} \sum_{j\neq i} p_{ij}(y+v)\times\\
&&\int_0^{\infty} \varphi(t+v,x,j,0)\alpha(x;s,i,v)\frac{\beta}{\sigma(i)\sqrt{v} s} dx dv,
\end{eqnarray}
\begin{eqnarray}\label{5ae}
\no \frac{\partial^2} {\partial s^2} \varphi(t, s, i, y)&=& \frac{ 1- F(T-t+y\mid i)}{1-F(y\mid i)}\frac{\partial^2 \eta_{i}(t,s)}{\partial s^2}+\int_0^{T-t} e^{-r(i)v}\frac{f(y+v\mid i)} {1-F(y\mid i)} \sum_{j\neq i} p_{ij} (y+v)\times\\
&&\int_0^{\infty} \varphi(t+v,x,j,0)\alpha(x;s,i,v) \frac{1}{s^2}\left(\frac{\beta^2}{\sigma^2(i) v} - \frac{\beta}{\sigma(i)\sqrt{v}} -\frac{1}{\sigma^2(i) v}\right)dx dv.
\end{eqnarray}
\noi From equations \eqref{4ae} and \eqref{5ae}, we get
\begin{eqnarray}\label{7ae}
\no \lefteqn{r(i) s \frac{\partial \varphi}{\partial s}+ \frac{1}{2}\sigma^2(i) s^2 \frac{\partial^2 \varphi}{\partial s^2}}\\
&=&\no \frac{ 1- F(T-t+y\mid i)}{1-F(y\mid i)}\left(r(i)s\frac{\partial\eta_{i}(t,s)}{\partial s}+\frac{1}{2}\sigma^2(i)s^2 \frac{\partial^2\eta_{i}(t,s)}{\partial s^2}\right)+\int_0^{T-t}e^{-r(i)v}\frac{f(y+v\mid i)}{1-F(y|i)}\\
&&\sum_{j\neq i} p_{ij}(y+v)\int_0^\infty \varphi(t+v,x,j,0)\alpha(x;s,i,v)\left(\frac{r(i)\beta}{\sigma(i)\sqrt{v}} +\frac{\beta^2}{2v}-\frac{\sigma(i)}{2\sqrt{v}}\beta-\frac{1}{2v}\right)dx dv.
\end{eqnarray}
\noi Finally, from equations \eqref{34}, \eqref{eq2}, \eqref{10ae}, \eqref{6ae} and \eqref{7ae} we get
\begin{eqnarray*}\label{8ae}
\lefteqn{\no\frac{\partial}{\partial t} \varphi(t, s, i, y)+\frac{\partial} {\partial y} \varphi(t, s, i, y)+r(i) s \frac{\partial }{\partial s}\varphi(t, s, i, y)+ \frac{1}{2}\sigma(i)^2(i) s^2\frac{\partial^2 }{\partial s^2}\varphi(t, s, i, y)}\\
&=&\no \frac{ 1- F(T-t+y\mid i)}{1-F(y\mid i)}\left[\frac{\partial\eta_{i}(t,s)}{\partial t} + r(i)s\frac{\partial\eta_{i}(t,s)}{\partial s}+\frac{1}{2}\sigma^2(i)s^2 \frac{\partial^2\eta_{i}(t,s)}{\partial s^2}\right] -\frac{f(y\mid i)}{1-F(y\mid i)}\times \\
&&\sum_{j\neq i} p_{ij}(y)(\varphi(t,s,j,0)-\varphi(t,s,i,y)) + r(i) \left(\varphi(t,s,i,y)- \frac{ 1- F(T-t+y\mid i)}{1-F(y\mid i)} \eta_{i}(t,s)\right)\\
&=& -\frac{f(y\mid i)}{1-F(y\mid i)}\sum_{j\neq i} p_{ij}(y)(\varphi(t,s,j,0)-\varphi(t,s,i,y))+r(i)\varphi(t,s,i,y).
\end{eqnarray*}
\noi Thus equation \eqref{p1} holds. \qed

\noi From Lemma \ref{lm2} and Proposition \ref{theo3} it follows that \eqref{p1}has a classical solution. We prove uniqueness in the following section.
\section{Uniqueness}
\begin{prop}\label{theo4} Assume (A1)(i)-(iii). Let $\varphi$ be a classical solution of \eqref{p1}. Then (i) $\varphi$ solves the integral equation \eqref{34}-\eqref{35}; (ii) $(s-K)^+\le \varphi(t,s,i,y) \le s$.
\end{prop}
\proof (i) Let $(\tilde\Omega,\tilde{\mathcal{F}}, \tilde P)$ be a probability space which holds a standard Brownian motion $\tilde W$ and the Poisson random measure $\tilde \wp$ independent of $\tilde W$. Let $\tilde{S}_t$ be the strong solution of the following SDE
\begin{eqnarray*}
d\tilde{S}_t = \tilde{S}_t(r(X_t)dt + \si(X_t)d\tilde W_t),\tab S_0>0
\end{eqnarray*}
where $X_t$ is the solution of \eqref{1} defined on $(\tilde\Omega,\tilde{\mathcal{F}}, \tilde P)$ and driven by $\tilde \wp$. Let $\tilde{ \mathcal{F}}_t$ be the underlying filtration generated by $(\tilde S_t, X_t)$ satisfying the usual hypothesis. We observe that the process $\{(\tilde{S}_t,X_t,Y_t)\}_t$ is Markov with infinitesimal generator $\mathcal{A}$, where
\begin{eqnarray*}
\mathcal{A} \varphi(s,i,y)= \frac{\partial \vf}{\partial y}  (s, i, y) + r(i) s \frac{\partial \vf}{\partial s} (s, i, y)+ \frac{1}{2} \si^2(i) s^2 \frac{\partial^2 \vf} {\partial
s^2} (s, i, y) + \sum_{j\neq i}\lambda_{ij}(y) \big(\varphi(s,j,0) -\varphi(s,i,y)\big)
\end{eqnarray*}
for every function $\varphi$ which is compactly supported $C^2$ in $s$ and $C^1$ in $y$. If $\vf$ is the classical solution of \eqref{p1} then by using the It\^{o}'s formula on $N_t := e^{-\int_0^t r(X_u)du} \varphi(t,\tilde{S}_t,X_t,Y_t)$, we get
\begin{eqnarray*}
  dN_t &=& e^{-\int_0^t r(X_u)du}\left(- r(X_t) \varphi (t,\tilde{S}_t,X_t,Y_t)+ \frac{\partial \varphi}{\partial t }(t,\tilde{S}_t,X_t,Y_t) + \mathcal{A} \varphi (t,\tilde{S}_t,X_t,Y_t)\right)dt + dM_t
\end{eqnarray*}
where $M_t$ is a local martingale. Thus from \eqref{p1} and above expression, $N_t$ is also an $ \tilde{ \mathcal{F}}_t$ local martingale. The definition of $N_t$ suggests that there are constants $k_1$ and $k_2$ such that $ |N_t|\le k_1 +k_2 \tilde S_t$ for each $t$, since $\vf$ has at most linear growth. Again, since the following expression
$$ \tilde S_t = \tilde S_0 \exp\left(\int_0^t (r(X_u)-\frac{1}{2}\sigma(X_u)^2)du + \int_0^t \sigma(X_u)d\tilde W_u\right)
$$
holds, one concludes that $\tilde S_t$ is a submartingale with finite expectation. Therefore Doob's inequality can be used to obtain $ \tilde E\sup_{s\in [0,t]}|N_s| < \infty$ for each $t$ where $\tilde E$ is the expectation w.r.t. $\tilde P$. Thus $\{N_t\}_t$ is a martingale. Hence
\begin{equation}\label{9}
\vf(t,\tilde{S}_t,X_t,Y_t)= e^{\int_0^t r(X_u)du}N_t = \tilde E[e^{\int_0^t r(X_u)du}N_T\mid \tilde{\mathcal{F}}_t]= \tilde E[e^{-\int_t^T r(X_u)du} (\tilde{S}_T-K)^+\mid \tilde{S}_t,X_t,Y_t].
\end{equation}
By conditioning at transition times and using the conditional lognormal distribution of $\tilde{S}_t$, we get
\begin{eqnarray*}
\lefteqn{\varphi(t,\tilde{S}_t,X_t,Y_t)}\\
&=& \tilde E[\tilde E[e^{-\int_t^T r(X_u)du} (\tilde{S}_T-K)^+ \mid \tilde{S}_t, X_t, Y_t, T_{n(t)+1}]\mid \tilde{S}_t, X_t, Y_t]\\
&=& P(T_{n(t)+1} > T\mid X_t,Y_t) \tilde E[e^{-\int_t^T r(X_u)du}  (\tilde{S}_T-K)^+\mid \tilde{S}_t, X_t, Y_t,T_{n(t)+1} > T] \\
&& + \int_0^{T-t} \tilde E[e^{-\int_t^T r(X_u)du} (\tilde{S}_T-K)^+\mid \tilde{S}_t, X_t, Y_t, T_{n(t)+1}=t + v] \frac{f(t-T_{n(t)}+v\mid X_t)} {1-F(Y_t\mid X_t)}dv\\
&=& \frac{ 1- F(T-T_{n(t)}\mid X_t)}{1-F(Y_t\mid X_t)} \eta_{X_t}(t,\tilde{S}_t)+ \int_0^{T-t} e^{-r(X_t)v} \frac{f(Y_t+v\mid X_t)} {1-F(Y_t\mid X_t)} \times\\ && \sum_j p_{X_tj}(Y_t +v) \int_0^{\infty} \tilde E[e^{-\int_{t+v}^T r(X_u)du} (\tilde{S}_T-K)^+\mid \tilde{S}_{t+v}=x,Y_{t+v}=0,\\
&& X_{t+v}=j, T_{n(t)+1}=t+v] \frac{e^{\frac{-1}{2}((\ln(\frac{x}{\tilde{S}_t})-(r(X_t) -\frac{\si^2(X_t)}{2})v) \frac{1}{\si(X_t)\sqrt v})^2}}{\sqrt{2\pi}\si(X_t)\sqrt{v} x} dx dv\\
&=& \frac{ 1- F(T-t+Y_t\mid X_t)}{1-F(Y_t\mid X_t)} \eta_{X_t}(t,\tilde{S}_t)+\int_0^{T-t} e^{-r(X_t)v} \frac{f(Y_t+v\mid X_t)} {1-F(Y_t\mid X_t)} \times\\
&& \sum_j p_{X_t j} (Y_t+v) \int_0^{\infty} \varphi(t+v,x,j,0) \frac{e^{\frac{-1}{2}((\ln(\frac{x}{\tilde{S}_t})-(r(X_t) -\frac{\si^2(X_t)}{2})v) \frac{1}{\si(X_t)\sqrt v})^2}}{\sqrt{2\pi} x\si(X_t)\sqrt v} dx dv.
\end{eqnarray*}
Finally by using irreducibility condition (A1), we can replace $( \tilde{S}_t, X_t,Y_t)$ by generic variable $(s,i,y)$ in the above relation and thus conclude that $\varphi$ is a solution of \eqref{34}-\eqref{35}. Hence the proof.

\noi (ii) We note that
$$\tilde E[e^{-\int_t^T r(X_u)du} (\tilde{S}_T-K)\mid \tilde{\mathcal{F}}_t] \le \tilde E[e^{-\int_t^T r(X_u)du} (\tilde{S}_T-K)^+\mid \tilde{\mathcal{F}}_t]\le \tilde E[e^{-\int_t^T r(X_u)du} \tilde{S}_T\mid \tilde{\mathcal{F}}_t].
$$
Now using martingale property of $e^{-\int_0^t r(X_u)du} \tilde{S}_t$, from \eqref{9} and above we get $\tilde{S}_t-K\le\vf(t,\tilde{S}_t,X_t,Y_t)\le \tilde{S}_t$. Again, from \eqref{9}, we know that $\vf$ is an expectation of a nonnegative quantity, hence nonnegative. Thus (ii) holds.\qed
\begin{theo}\label{theo1}
The initial-boundary value problem \eqref{p1}has a unique classical solution in the class of functions with at most linear growth.
\end{theo}
\proof Existence follows from Lemma \ref{lm2} and Proposition \ref{theo3}. For uniqueness, first assume that $\vf_1$ and $\vf_2$ are two classical solutions of \eqref{p1}in the prescribed class. Then using Proposition \ref{theo4}, we know that both also solve \eqref{34}-\eqref{35}. But from Lemma \ref{lm2}, there is only one such in the prescribed class. Hence $\vf_1=\vf_2$.\qed
\begin{rem} The above theorem can also be proved in a different manner which heavily depends on the mild solution techniques \cite{PA} and Proposition 3.1.2 of \cite{CFMW}. Such an alternative approach is taken in \cite{AGMKG} to establish well-posedness of a special case of \eqref{p1}. The reason for adopting the present approach is that, it enables us to establish the equivalence between the PDE and an IE in one go. This in turn suggests an alternative expression of partial derivative of the solution. In the next section the importance of such representation is explained.
\end{rem}

\section{Pricing and optimal hedging}
In this section we consider European call option on the stock dynamics as given in Section 2. The terminal claim, associated to that option is $(S_T-K)^+$. We show that its locally risk minimizing price at time $t(\le T)$ can be written in terms of the solution of \eqref{p1}. We further show that the corresponding optimal hedging has an integral representation in terms of the price function.

\begin{theo}\label{theo5} Let $\vf$ be the unique classical solution of \eqref{p1}in the class of functions with at most linear growth.
\begin{enumerate}
  \item Let $(\xi,\varepsilon)$ be given by
  \begin{equation}\label{VI3.20}
 \xi_t :=\frac{\partial\varphi(t,S_t,X_{t-},Y_{t-})}{\partial s} \txt{ and } \varepsilon_{t} := e^{-\int_{0}^{t}r(X_{u})du} (\varphi(t,S_t,X_{t},Y_{t})-\xi_{t}S_{t}).
 \end{equation}
Then $(\xi,\varepsilon)$ is the optimal admissible strategy.
  \item $\varphi(t,S_t,X_t,Y_t)$ is the locally risk minimizing price of $(S_T-K)^+$ at time $0\le t\le T$.
\end{enumerate}
\end{theo}

\proof Under the market model, the mean variance tradeoff (MVT) process $\hat{K}_t$ (as defined in Pham et al \cite{PH}) takes the following form
\begin{equation*}
\hat{K}_t=\int_0^t\left(\frac{\mu(X_s)-r(X_s)}{\sigma(X_s)}\right)^2 ds.
\end{equation*}
Hence $\hat{K}_t$ is bounded and continuous on $[0,T]$. We also know that $S_t$ has almost sure continuous paths. Since, $H^*\in L^2(\Omega, {\cal F}, P)$ for $H=(S_T-K)^+$ we apply corollary 5 and Lemma 6 of \cite{PH} to conclude that $H^*$ admits F\"{o}llmer-Schweizer decomposition \eqref{eq1} with an integrand $\xi^{H^*}$ satisfying A2(i) and $L^{H^*}$ being square integrable. Therefore, to prove the theorem it is sufficient to show that
\begin{itemize}
  \item[(a)] $\varphi(t,S_t,X_{t},Y_{t})= \varepsilon_{t} B_t + \xi_{t}S_{t}$ for all $t\le T$;
  \item[(b)] $P(\varphi(t,S_t,X_{t},Y_{t}) \ge 0 \forall t\le T)=1$
  \item[(c)] $\frac{1}{B_t} \varphi(t,S_t,X_{t-},Y_{t-})= H_0+\int^{t}_{0}{\xi_t}dS^*_t+L_t$ for all $t\le T$;
  \item[(d)] there exists $\mathcal{F}_0$ measurable $H_0$ and $\mathcal{F}_T$ measurable $L_T$ such that $L_t:= E[L_T\mid \mathcal{F}_t]$ is orthogonal to $\int_0^t \sig(X_t)S^*_t dW_t$ i.e., the martingale part of $S^*_t$ and $H^*=H_0+\int^{T}_{0}{\xi_t}dS^*_t+L_T$;
\end{itemize}
where $\vf$ is the unique classical solution of \eqref{p1}  in the prescribed class and $(\xi,\varepsilon)$ is as in \eqref{VI3.20}.

\noi  From the definition of $\varepsilon_t$ in \eqref{VI3.20}, (a) follows. In Lemma \ref{lm2} and Proposition \ref{theo4}, it is shown that $\vf$ is a non-negative function. Hence (b) holds. Next we show the condition (c) and (d). We apply It\^{o}'s formula to $e^{-\int_{0}^{t}r(X_{u})du} \varphi(t,S_{t},X_{t},Y_{t})$ under the
measure $P$ and use \eqref{8}, \eqref{p1}  and \eqref{1}
to obtain after suitable rearrangement of terms,  for all $t<T$
\begin{eqnarray*}\label{VI3.24}
e^{-\int_{0}^{t}r(X_{u})du} \varphi(t,S_{t},X_{t},Y_{t}) & = &
\varphi(0,S_{0},X_{0},Y_{0})+\int_{0}^{t}\frac{\partial
\varphi(u,S_{u},X_{u-},Y_{u-})}{\partial s} d S^*_{u}\no
+\int_{0}^{t}e^{-\int_{0}^{u}r(X_{v})dv}\\
& &\hspace{-5mm}\int_{\mathbb{R}} [\varphi(u,S_{u},X_{u-} +h(X_{u-},Y_{u-},z),Y_{u-}-g(X_{u-},Y_{u-},z))\\
&&-\varphi(u,S_{u},X_{u-},Y_{u-})]{\hat{\wp}}(du,dz)
\end{eqnarray*}
where the last integral is w.r.t. the compensator of $\wp$. We set
\begin{eqnarray*}
L_t&:=& \int_{0}^{t}e^{-\int_{0}^{u}r(X_{v})dv}\int_{\mathbb{R}} [\varphi(u,S_{u},X_{u-} +h(X_{u-},Y_{u-},z),Y_{u-}-g(X_{u-},Y_{u-},z)) \\
&&-\varphi(u,S_{u},X_{u-},Y_{u-})]{\hat{\wp}}(du,dz).
\end{eqnarray*}
We note that Proposition \ref{theo4} (ii) implies that the integrand in above expression is bounded by $K$. Therefore, $L_t$ being an integral w.r.t. a compensated Poisson random measure, is a martingale. Again the independence of $W_t$ and $\wp$ implies the orthogonality of $L_t$ to the martingale part of $S^*_t$. Thus, we obtain the following F\"{o}llmer-Schweizer  decomposition by letting $t \uparrow T $,
\begin{equation}\label{VI3.25}
B_T^{-1}(S_{T}-K)^+ = \varphi(0,S_{0},X_{0},Y_{0}) + \int^{T}_{0}{\xi_t}dS^*_t + L_T.
\end{equation}
Thus (c) and (d) hold.\qed
\begin{theo}\label{theo6}
Let $\varphi$ be the unique solution of \eqref{p1}. Set
\begin{eqnarray}\label{5}
\no\psi(t,s,i,y) &:=& \frac{ 1- F(T-t+y\mid i)}{1-F(y\mid i)} \frac{\partial\eta_{i}(t,s)}{\partial s}+\int_0^{T-t}e^{-r(i)v} \frac{f(y+v\mid i)} {1-F(y\mid i)} \times\sum_{j\neq i} p_{ij}(y+v) \\
\no&& \int_0^{\infty} \varphi(t+v,x,j,0) \frac{e^{\frac{-1}{2}((\ln(\frac{x}{s})-(r(i) -\frac{\si^2(i)}{2})v) \frac{1}{\si(i)\sqrt v})^2}}{\sqrt{2\pi} xs\sigma(i)\sqrt{v}}\frac{\left(\ln(\frac{x}{s})-(r(i) -\frac{\si^2(i)}{2})v\right)}{\sigma(i)^2v} dx dv\\
\end{eqnarray}
where $(t,s,i,y)\in \mathcal{D}$.
Then $\psi(t,s,i,y) =\frac{\pa}{\pa s}\vf(t,s,i,y)$.
\end{theo}

\proof We need to show that $\psi$ (as in \eqref{5}) is equal to $\frac{\partial\varphi}{\partial s}$. Indeed, one obtains the RHS of \eqref{5} by differentiating the right side of \eqref{34} with respect to $s$. Hence the proof. \qed

\begin{remark}\label{rem1}
It is well known that in a numerical differentiation, an isolated perturbation gets enhanced whereas in a numerical integration the same gets reduced. In \eqref{5}, the function $\psi$, a partial derivative of $\varphi$, is given by an integration involving $\varphi$. Thus the above theorem essentially provides a robust way to find $\frac{\pa}{\pa s}\vf(t,s,i,y)$ using a step by step quadrature method. Further we show that $\frac{\pa\vf}{\pa s}$ constitutes the locally risk minimizing hedging of a contingent claim whose price is given by $\vf$. We would like to mention that the price of other types of options such as barrier options, composite options, basket options etc can also be represented similarly by imposing appropriate terminal and boundary conditions. Thus hedging strategy of those would also have similar integral representations. Needless to mention that this finding holds true for Markov modulated market model, a special case of the present model.
\end{remark}
\section{Conclusion}
To our knowledge, the question of robust computation of hedging for regime switching market is not addressed in the literature. In this paper an integral representation of optimal hedging is obtained which leads to a robust scheme of computation of hedging. Besides, there is another important contribution. The integral equation established here, exhibits how option price in regime switching case depends on that in Black-Scholes counterpart. In equation \eqref{5}, we indeed express the optimal hedging using Black-Scholes delta hedging and some additional terms. It is also important to note that $\lambda_{ij}$ are the only additional(functional) parameter from the Black-Scholes case. Using some nonparametric estimation procedure, as in \cite{OL}, it is possible to have a consistent estimate of $\lambda_{ij}$ functions, if the regimes are observed. In a recent paper \cite{GoN} the convergence of corresponding approximate option price functions is investigated.


\begin{thebibliography}{99}

\bibitem{CFMW} Arendt W., Batty C., Hieber, M. and Neubrander, F., Vector-valued Laplace Transforms and Cauchy Problems, Birkhauser 2001.

\bibitem{BAS} Basak G. K., Ghosh Mrinal K. and Goswami A., Risk minimizing option pricing for a class of exotic options in a Markov-modulated market,  Stoch. Ann. App. 29:2(2011), 259-281.

\bibitem{BS}Black F. and Scholes M., The pricing of options and corporate liability, Journal of Political Economy, Vol. 81, No. 3 (1973), 637-659.

\bibitem{BE}Buffington J. and Elliott R. J., American options with regime switching, Intl. J. Theor. Appl. Finance 5(2002), 497-514.

\bibitem{DES} Deshpande A. and Ghosh M. K., Risk minimizing option pricing in a regime switching market, Stoch. Ann. App. 26(2008).

\bibitem{DKR} DiMasi G. B., Kabanov Y. and Runggaldier W. J., Mean-Variance hedging of options on stocks with Markov volatility. Theory Probab. Appl., Vol. 39 (1994), 173-181.

\bibitem{ECS} Elliott R.J., Chan L. and Siu T.K., Option pricing and Esscher transform under regime switching, Annals of Finance 1,
423-432 (2005).

\bibitem{FS} F\"{o}llmer H. and Schweizer, M., \emph{Hedging of Contingent Claims under Incomplete Information}, Applied Stochastic Analysis, Stochastics Monographs, vol. 5 (1991), 389-414. 

\bibitem{AGMKG} Ghosh M. K. and Goswami A., Risk minimizing option pricing in a semi-Markov modulated market, SIAM J. Control Optim. 48(2009), 1519-1541.

\bibitem{GhS} Ghosh M. K. and Saha, S., Stochastic processes with age-dependent transition rates, Stoch. Ann. App. 29(2011), 511-522.

\bibitem{GoN} Goswami, A. and Nandan, S., Convergence of estimated option price in a regime switching market, Indian J. Pure Appl. Math. to appear.

\bibitem{GoS} Goswami, A. and Saini, R. K., Volterra equation for pricing and hedging in a regime switching market, Cogent Economics and Finance 2(2014), no.1, 1-11.

\bibitem{GZ}Guo X. and Zhang Q., Closed form solutions for perpetual American put options with regime switching, SIAM J. Appl. Math 39(2004), 173-181.

\bibitem{HE} Heston, Steven L., A closed-form solution for options with stochastic volatility with applications to bond and currency options, Review of financial studies 6.2 (1993): 327-343.

\bibitem{JH} Hunt J. and Devolder P., Semi-Markov regime switching interest rate models and minimal entropy measure, Physica A: Statistical Mechanics and its Applications 390, 15(2011), 3767-3781.

\bibitem{RR} Joberts A. and Rogers L. C. G., Option pricing with Markov-modulated dynamics, SIAM J. Control Optim. 44(2006), 2063-2078.

\bibitem{KK} Kallianpur Gopinath and Karandikar Rajeeva L., Introduction to Option Pricing Theory, Birkhauser Boston, 2000.

\bibitem{MR} Mamon R. S. and Rodrigo M. R., Explicit solutions to European options in a regime switching economy, Operations Research Letters 33(2005), 581-586.

\bibitem{OL}
Ouhbi, B. Limnios, N., Nonparametric estimation for semi-Markov processes based on its hazard rate functions. Stat. Inference Stoch. Process., 2(1999), 151-173.

\bibitem{PA} Pazy A., Semigroups of Linear Operators and Applications to Partial Differential Equations, Springer-Verlag,
1983.

\bibitem{PH} Pham H., Rheinl\"{a}nder T. and Schweizer, M., Mean-variance hedging for continuous processes: new proofs and examples, Finance Stoch.(1998) 173-198.

\bibitem{S2} Schweizer M., A Guided Tour through Quadratic Hedging Approaches, E. Jouini, J. Cvitani\'{c}, M. Musiela (eds.), Option Priing Interest Rates and Risk Management, Cambridge University Press (2001), 538-574.
\end{thebibliography}
\end{document}